\documentclass{article}
\usepackage[english]{babel}
\usepackage[a4paper,top=2cm,bottom=2cm,left=3cm,right=3cm,marginparwidth=1.75cm]{geometry}
\usepackage{tabularx}
\usepackage{multirow}
\usepackage{setspace}
\usepackage{amssymb}
\usepackage{amsmath}
\usepackage{siunitx}
\PassOptionsToPackage{hyphens}{url}\usepackage{hyperref}
\usepackage[utf8]{inputenc}
\usepackage[right]{lineno}
\usepackage{csquotes}
\usepackage{booktabs}
\usepackage{longtable}
\usepackage{adjustbox}
\usepackage{array}
\usepackage{url}
\usepackage{titlesec}
\usepackage{authblk}
\usepackage{xcolor} 

\usepackage{indentfirst}
\usepackage{dsfont}

\usepackage[english]{babel}

\author[1]{Diego R. Rivera}
\author[1]{Ernesto Castillo}
\author[2]{Felipe Galarce\footnote{\normalsize{\textbf{Corresponding author: \url{felipe.galarce@pucv.cl}}}}}
\author[3]{Douglas R.Q. Pacheco}

\affil[1]{\small Departamento de Ingeniería Mecánica, Universidad de Santiago, Santiago, Chile.}
\affil[2]{\small School of Civil Engineering, Pontificia Universidad Católica de Valparaíso, Valparaíso, Chile}
\affil[3]{\small Chair for Computational Analysis of Technical Systems, RWTH Aachen University, Germany.}

\title{Improving performance estimation of a PCM-integrated solar chimney through reduced-order based data assimilation}

\begin{document}
\maketitle

\begin{abstract}    
\textbf{Purpose} – This study assesses a data-assimilation framework based on reduced-order modeling (ROM–DA), complemented with a data-filling strategy, to reconstruct the dynamic temperature fields of a PCM-integrated solar chimney from scarce temperature measurements, to improve the estimation accuracy of outlet airflow velocity. 

\textbf{Design/Methodology/Approach} – A regularized least-squares formulation is employed to estimate temperature fields within an inclined solar chimney configuration using RT-42 as the phase-change-material (PCM). The method combines: (i) a reduced-order model derived from high-fidelity finite-volume simulations of unsteady conjugate heat transfer with liquid–solid phase change and surface radiation, and (ii) three measurement datasets of 22, 135, and 203 points. The datasets are expanded using a hybrid data-filling strategy (boundary-layer and bi-cubic interpolations). The reconstructed temperature fields are subsequently assimilated into the thermally-coupled forward solver to enhance the airflow velocity prediction.

\textbf{Findings} – The ROM–DA framework accurately reconstructed dynamic temperature fields in both the air and PCM domains using synthetic measurements with relative errors below 10\% and 3\% for the initial and expanded sensor sets, respectively. When applied to real measurements, the framework improved the fidelity of the local temperature evolution in both domains relative to the forward model. Increasing the number of sensors did not significantly improve local temperature accuracy but enhanced the dynamic estimation of the local outlet velocity by reducing the RMS error by 20\%.

\textbf{Originality/Value} – This is the first application of a ROM–DA framework to a coupled multiphysics solar chimney with PCM integration. The study also assesses hybrid data-enhancement strategies that improve measurement quality. 

\end{abstract}


\section{Introduction} \label{sec:introduction} 

Heating, ventilation, and air conditioning (HVAC) systems have seen a steady rise in global demand over recent decades, driven by demographic and economic growth in emerging economies, where higher living standards have increased the need for indoor thermal comfort as a key factor for health, productivity, and well-being (\cite{iea_2025_stayingcool}). Concurrently, the growing frequency and intensity of temperature extremes associated with climate change—projected to persist at least until 2050—(\cite{falchetta_2024_natcomm_cooling}) are reinforcing the world’s dependence on HVAC technologies.

This widespread reliance on conventional HVAC systems carries a heavy environmental and energy cost. They account for nearly 38\% of total building energy use, 12\% of global end-use energy demand, and roughly a quarter of worldwide energy-related CO$_2$ emissions (\cite{gonzalez-torres_2022_energyreports}). These values reveal the urgent need for efficient and sustainable alternatives. Passive HVAC strategies that exploit renewable resources, such as solar-assisted designs and thermal-mass optimization, can reduce or even replace active mechanical systems, lowering building heating and cooling loads by up to 40\% (\cite{cao_2016_enbuild}) while supporting global decarbonization targets.

Among passive solar HVAC technologies, Trombe walls and solar chimneys are well-established designs. Both consist of a collector with a glass cover, an air channel, and an absorptive plate, where absorbed solar energy drives buoyancy-induced airflow that provides heating or ventilation without mechanical aid. However, natural variations in solar irradiance due to clouds or diurnal cycles limit their thermal and flow stability. A common enhancement is the integration of high-thermal-inertia media such as phase change materials (PCMs), which store and release latent heat through quasi-isothermal transitions, thereby damping temperature and flow oscillations. The use of PCMs in solar-assisted HVAC systems has been extensively investigated over the past four decades (\cite{vargas-lopez_2019_energy, omara_2021_solener_review, rashid_2025_jenergystorage_review}), showing higher day-night air exchange rates, improved thermal and exergy efficiencies, and smoother indoor temperature profiles. Nevertheless, identifying suitable PCMs and optimizing system design under variable climatic and operating conditions remain challenging, while purely experimental approaches are often costly and time-consuming.

Numerical modeling has therefore become essential for studying solar-driven HVAC systems with PCMs. Their multiphysics nature—buoyancy-driven convection, conjugate heat transfer across solid-air-PCM interfaces, surface radiation, and transient melting and solidification—requires high-fidelity computational fluid dynamics (CFD). Reliable predictions demand fine spatial and temporal resolution to capture the strong nonlinear coupling between momentum and energy transport. Several studies (\cite{tiji_2020_jobe, qasim_2024_ijft, salari_2020_solener, safari_2014_ees, kong_2020_solener}) have shown that such high-resolution CFD frameworks can accurately reproduce convective heat transfer and phase-change front dynamics. However, those simulations are computationally intensive, particularly for diurnal-cycle analyses, and their accuracy is constrained by uncertainties in thermophysical properties, boundary and initial conditions, as well as by discretization-induced bias and numerical diffusion errors (\cite{jimenez-xaman_2019_enbuild}).

In this context, data assimilation (DA) has emerged as a powerful strategy to systematically integrate our physical knowledge of the phenomena, in disregard of the parameter uncertainty and modeling limitations, with experimental data in complex thermo-fluid systems. By merging the governing equations with measurements and statistical error models, DA transforms purely predictive solvers into adaptive estimators that continuously correct their state, thereby reducing systematic bias and posterior uncertainty (\cite{he_2025_experfluids, evensen2009data}). Ill-conditioning and computational-cost problems of high-fidelity solvers in DA can be alleviated with reduced-order models (ROMs) that preserve the underlying physics while drastically reducing the number of degrees of freedom (\cite{BCDDPW2017}). In combination, ROM-based DA enable near-real-time state estimation and control with fidelity sufficient for physics-aware digital twins (\cite{arcucci_2023_compfluids, eshaghi_2024_compstruct}).

Although DA has long been a cornerstone of weather forecasting and climate modeling (\cite{martin_2025_stateplanet}), its systematic application in engineering disciplines has expanded significantly in recent years. Notable advances include the assimilation of buoyancy-driven convection data (\cite{farhat_2018_jsci}), the integration of tomographic PIV measurements into high-fidelity simulations (\cite{bauer_2022_experfluids}), and recent HVAC-focused studies addressing mixing and natural ventilation flows (\cite{liu_2024_jtsea, wu_2018_fluids, qian_2025_buildenv_da_validation}). Collectively, these works demonstrate the capability of DA frameworks to reconstruct temperature and velocity fields from sparse or noisy measurements, improving physical consistency and predictive accuracy. However, their extension to solar-driven HVAC systems with PCM integration remains largely unexplored, particularly concerning the simultaneous reconstruction of airflow and PCM domains and the impact of measurement density or data-filling strategies on assimilation performance.

In this context, many critical knowledge gaps remain in the application of ROM--DA frameworks to engineering systems—particularly in solar-driven thermal configurations incorporating PCMs. It is still unclear how to reconstruct coupled states across multiple domains (airflow and PCM) or whether accurate airflow prediction can be achieved solely through temperature assimilation. Moreover, the role of measurement density and the potential of simple gap-filling techniques to enhance assimilation robustness have not been systematically examined. 
To address these issues, this study evaluates a ROM–DA framework for a solar collector composed of glass, air, and PCM layers, supported by a high-fidelity CFD solver validated against experimental data, with the objective of improving ventilation-rate prediction. We also investigate how physics- and data-driven filling techniques influence the stability, reliability, and computational efficiency of the assimilation process. To the best of our knowledge, this is the first work to integrate a ROM–DA framework into a PCM-enhanced solar chimney, bridging advanced data-driven techniques with multiphysics building energy systems.

The manuscript is organized as follows. Section \ref{sec:Method} states the methodology, including the physical configurations, governing equations, and initial/boundary conditions. Section \ref{sec:Numerical} summarizes the numerical methods, encompassing the forward solver and data assimilation frameworks. Section \ref{sec:Results} presents and discusses the main results. Finally, Section \ref{sec:Conclusions} draws concluding remarks.

\section{Methodology}\label{sec:Method} 

\subsection{Physical configuration} 

\begin{figure}[h]
    \centering
    \includegraphics[width=1\linewidth]{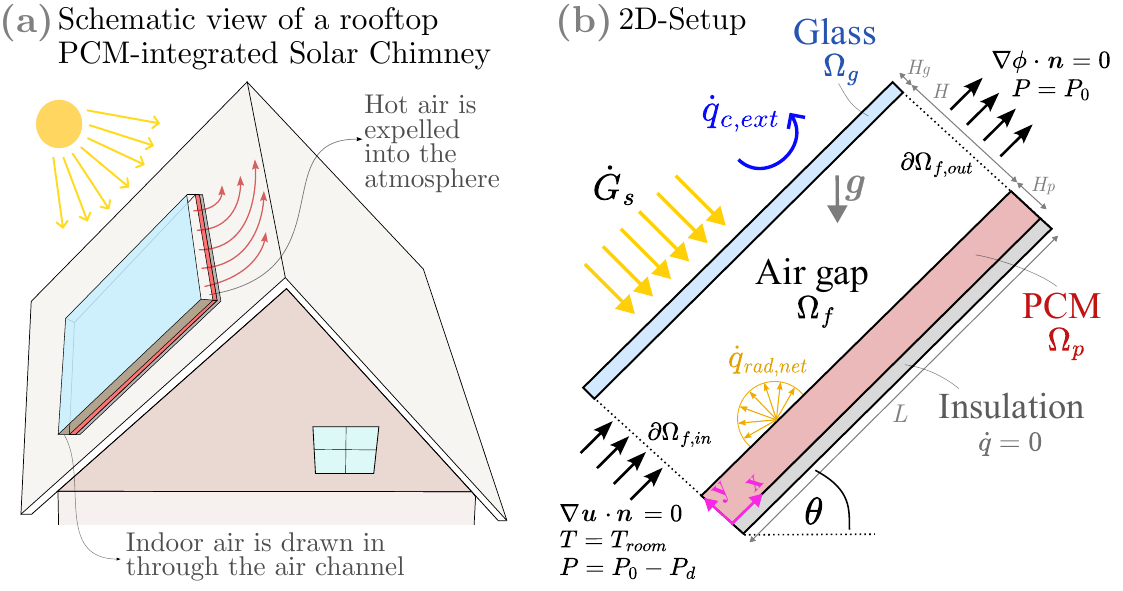}
    \caption{(a) Schematic representation of an inclined solar chimney. (b) Two-dimensional setup.}
    \label{fig:physical_situation}
\end{figure}

The solar chimney we analyzed corresponds to a flat-plate collector system installed at an inclination angle on the rooftop of buildings, as illustrated in Fig.~\ref{fig:physical_situation}(a). Following common modeling practice, the analysis is restricted to a two-dimensional domain at the mid-plane, as presented in Fig.~\ref{fig:physical_situation}(b), which adequately captures the dominant physical processes \cite{nguyen_2024_jbuildphycs}.  

The configuration replicates the experimental setup reported by \cite{HUANG2024130154}, consisting of three flat layers of total length $L=1$~m: a glass cover, an intermediate air gap of 0.3~m, and a 0.03~m thick aluminum absorber plate painted black and thermally insulated on its rear surface. The absorber plate also serves as a macro-encapsulation container for the commercial PCM RT42, whose thermophysical properties are summarized in Table~\ref{tab:pcm_props}.  

The experimental setup comprises a charging stage under constant solar irradiation $\dot{G}_{S}$ ranging from 100 to 600~W/m$^2$, emulated by halogen lamps, followed by a discharge stage monitored until the system reaches the ambient temperature $T_{\mathrm{amb}}=22\,^{\circ}$C. Tests were conducted at three inclination angles: $\theta=30^{\circ}$, $45^{\circ}$, and $60^{\circ}$. Since the original reference does not provide detailed boundary specifications, physically consistent assumptions are adopted based on the reported inlet temperature evolution. The collector is considered to be inside a room, while the glass cover exchanges heat with the surrounding air solely through natural convection. These assumptions are consistent with typical indoor solar collector experiments and ensure a realistic thermal boundary treatment.

\begin{table}[h!]
\centering
\begin{tabular}{l |c| c c c c}
    \hline
    \textbf{Property} & \textbf{Value} & \multicolumn{4}{c}{$C_p(T) = A T^4 + B T^3 + C T^2 + D T + E$} \\
    \hline
    Phase-change range ($T_{\mathrm{sol}}/T_{\mathrm{liq}}$), $^\circ$C & $37 / 43$ & & $T < T_1$ & $T_1 \leq T \leq T_2$ & $T > T_2$ \\
    Latent heat ($H_{\mathrm{pc}}$), kJ/kg & $165$ & A & $5.22\times 10^{-2}$ & $-2.46\times 10^{0}$ & $1.00\times 10^{0}$ \\
    Density ($\rho_{\mathrm{sol}}/\rho_{\mathrm{liq}}$), kg/m$^3$ & $880 / 760$ & B & $-2.46\times 10^{-1}$ & $6.29\times 10^{1}$ & $-4.14\times 10^{1}$ \\
    Thermal conductivity, W/m$^\circ$C & $0.20$ & C & $6.54\times 10^{-1}$ & $-6.00\times 10^{2}$ & $6.43\times 10^{2}$ \\
    Viscosity, Pa·s & $0.0235$ & D & $-2.86\times 10^{-1}$ & $2.52\times 10^{3}$ & $-4.44\times 10^{3}$ \\
    Thermal expansion coeff., K$^{-1}$ & $0.0005$ & E & $3.80\times 10^{0}$ & $-3.94\times 10^{3}$ & $1.15\times 10^{4}$ \\
    \hline
\end{tabular}
\caption{Thermophysical properties of the commercial PCM RT42. The effective specific heat extracted from the enthalpy-temperature curves is fitted to a fourth-degree piecewise polynomial function, with $T_1=$ 39 $^{\circ}$C and $T_2=$ 32.5 $^{\circ}$C.}
\label{tab:pcm_props}
\end{table}


\subsection{Governing equations} 

The mathematical model of the glass-air-PCM system is defined separately for each material domain. In the air and PCM regions, the flow is governed by the incompressible Navier–Stokes equations coupled with the energy equation through the Boussinesq approximation to account for buoyancy effects. For compactness, the equations for each domain are expressed using a subscript $i$, which identifies the material (air or PCM) and acts as a switch to activate or deactivate the terms associated with each region:
\begin{equation}
  \begin{aligned}
        \nabla \cdot \boldsymbol{u}_i &= 0, 
        && \text{in } \Omega_i \times (0,t_{\mathrm{end}}],\\
        \rho_i \partial_{t} \boldsymbol{u}_i + \rho_i (\boldsymbol{u}_i \cdot \nabla) \boldsymbol{u}_i 
        - \mu_i\nabla^2  \boldsymbol{u}_i + \nabla P_i &= \boldsymbol{f}_b + \boldsymbol{f}_{\mathrm{sol}},
        && \text{in } \Omega_i \times (0,t_{\mathrm{end}}],\\
        \rho_i C_i \partial_{t} T_i + \rho_i C_i \boldsymbol{u}_i \cdot \nabla T_i  
        - \lambda_i \nabla^2 T_i &= 0,
        && \text{in } \Omega_i \times (0,t_{\mathrm{end}}],
  \end{aligned}    
\end{equation}
where subscript $i$ denotes the material domain, with $i=f$ for the air and $i=p$ for the PCM. The state variables are the velocity vector $\boldsymbol{u}$, pressure $P$, and temperature $T$. Thermophysical properties $\rho$, $\mu$, $C$, and $\lambda$ (density, dynamic viscosity, specific heat, and thermal conductivity, respectively) are considered constant for the air and temperature-dependent for the PCM according to the enthalpy–porosity formulation.

Latent heat storage is represented through an effective temperature-dependent specific heat obtained from the enthalpy–temperature curves provided by the PCM-RT42 supplier's data sheet (see Table~\ref{tab:pcm_props}).
Density and thermal conductivity are linearly interpolated between the solid and liquid phases according to the liquid fraction $f_{\mathrm{pc}}$. 
The hydrodynamic transition between solid and liquid states is modeled through a porosity-dependent source term in the momentum equation $\boldsymbol{f}_{\mathrm{sol}}$. Both quantities are defined as:
\begin{equation}
\begin{aligned}
f_{\mathrm{pc}} &= 
\begin{cases} 
    0, & T_p < T_{\mathrm{sol}}, \\[6pt]
    \dfrac{T_p - T_{\mathrm{sol}}}{T_{\mathrm{liq}} - T_{\mathrm{sol}}}, & T_{\mathrm{sol}} \leq T_p \leq T_{\mathrm{liq}}, \\[10pt]
    1, & T_p > T_{\mathrm{liq}},
\end{cases}
& \qquad
\boldsymbol{f}_{\mathrm{sol}} = C_M \,\frac{(1-f_{\mathrm{pc}})^2}{f_{\mathrm{pc}}^3+\varepsilon}\,\boldsymbol{u},
\end{aligned}
\label{phasechange_model}
\end{equation}
where $T_{\mathrm{sol}}$ and $T_{\mathrm{liq}}$ denote the solidus and liquidus temperatures of the PCM, $C_M=10^5$ is the mushy-zone constant, and $\varepsilon=10^{-4}$ prevents division by zero.  

Heat transfer within the glass is governed exclusively by the energy equation with constant thermophysical properties, as only conduction is relevant:
\begin{equation}
    \rho_g C_g \partial_{t} T_g  - \lambda_g \nabla^2 T_g + \partial_{y}\Theta_g = 0, 
    \qquad \text{in } \Omega_g \times (0,t_{\mathrm{end}}],
\end{equation}
where solar energy absorption is represented by an exponential attenuation function:
\begin{equation}
\Theta_g = \dot{G}_{S}\,\exp[-S_g(H_g-y)],
\end{equation}
with $S_g=20$~m$^{-1}$ being the extinction coefficient and $H_g=0.006$~m the glass thickness.  

Surface thermal radiation between the absorber plate and the glass cover is modeled using the net-radiation method. 
For each surface element $j$, the net radiative heat flux is computed as:
\begin{equation}
\dot{q}_{\mathrm{rad,net}}^{(j)} = \dot{q}_{\mathrm{rad,out}}^{(j)} - \dot{q}_{\mathrm{rad,in}}^{(j)},
\end{equation}
where the outgoing and incoming fluxes are:
\begin{equation}
\dot{q}_{\mathrm{rad,out}}^{(j)} = \epsilon\,\sigma T_{j}^{4} + (1-\epsilon)\,\dot{q}_{\mathrm{rad,in}}^{(j)},
\qquad
\dot{q}_{\mathrm{rad,in}}^{(j)} = \sum_{i=1}^{N_s} \dot{q}_{\mathrm{rad,out}}^{(i)}\,F_{ji},
\label{radiation_model}
\end{equation}
with $\epsilon$ being the emissivity of the surface, $\sigma$ the Stefan–Boltzmann constant, and $F_{ji}$ the view factors computed via the string method. The emissivities of the glass and black-painted absorber are 0.85 and 0.95, respectively, while the short-wave transmittance of the glass is $\tau_g=0.78$. 
Since the incoming and outgoing fluxes are interdependent, Eqs.~\eqref{radiation_model} are solved iteratively within each outer iteration of the coupled solver, with a stopping criterion for the inner iterations $k$ as
\[
\max\big| \dot{q}_{\mathrm{rad,net}}^{(k+1)} - \dot{q}_{\mathrm{rad,net}}^{(k)} \big| \leq 10^{-5}.
\]

\subsection{Initial and boundary conditions} 

Boundary conditions applied to the glass–air–PCM system are summarized in Table~\ref{tab:boundaryc}. 
At the inlet boundary $\partial\Omega_{\mathrm{f,in}}$, the velocity satisfies an open boundary condition, the temperature is prescribed at the room's temperature, and the pressure is defined as the sum of static and dynamic contributions, as typically used in buoyancy-driven flow configurations. 
At the outlet boundary $\partial\Omega_{\mathrm{f,out}}$, velocity follows a developed-flux condition corrected by a continuity factor $m_c=\sum \dot{m}_{\mathrm{in}}/\sum \dot{m}_{\mathrm{out}}$, the temperature obeys an outflow condition, and the pressure is fixed to the static reference value $P_0=1$~atm.
Finally, the conditions at the air-glass and air-PCM interfaces are defined by energy balances between heat fluxes by conduction, convection, and radiation; the convective coefficient at the outer wall of the glass $h_{\mathrm{ext}}$ is calculated with a correlation for laminar natural convection in an inclined flat plate (\cite{perovic_2016_thermsci}).

Based on the measured inlet temperature, we inferred that the experimental room was not maintained at constant temperature and exhibits a thermal mass inertia $I_r$ associated with the charging and discharging stages. 
Accordingly, the temporal evolution of the room's temperature is modeled by a global energy balance, in which the absorbed energy depends on the incident radiation $\dot{G}_S$, and the heat loss occurs by convective exchange with the ambient air:
\begin{equation}
    I_{r}\,\partial_tT_{\mathrm{room}} = \dot{G}_S - C_1\,(T_{\mathrm{room}}-T_{\mathrm{amb}}),
\end{equation}
where $I_r$ represents the thermal mass inertia and $C_1$ the convective exchange coefficient. 
Both parameters are adjusted to reproduce the experimental inlet temperature evolution during the charging stage, yielding $I_r=188.5$~kJ/m$^2$K and $C_1=55$~W/m$^2$K.

\begin{table}[h]
    \centering
        \begin{tabular}{c c c c} 
        \hline
         Boundary domain & Velocity & Temperature & Pressure \\
        \hline
         $\partial\Omega_{\mathrm{f,in}}$ & $\nabla\boldsymbol{u}_f\cdot \boldsymbol{n}=0$ & $T_f=T_{\mathrm{room}}$ & $P_f=P_0-\tfrac{1}{2}\rho U^2$ \\
        \hline
         $\partial\Omega_{\mathrm{f,out}}$ & $m_{c}\nabla\boldsymbol{u}_f\cdot \boldsymbol{n}=0$ & $\nabla T_f\cdot \boldsymbol{n}=0$ & $P_f=P_0$ \\ 
        \hline
         $\partial\Omega_{f-g}$ & $\boldsymbol{u}_f=0$ & $\lambda_f \nabla T_f \cdot \boldsymbol{n}= \lambda_{g}\nabla T_{g} \cdot \boldsymbol{n} + \dot{q}_{g,\mathrm{rad,net}}$ & — \\
        \hline
         $\partial\Omega_{f-p}$ & $\boldsymbol{u}_f=0,\, \boldsymbol{u}_p=0$ & $\lambda_f \nabla T_f \cdot \boldsymbol{n}= \lambda_{p}\nabla T_{p} \cdot \boldsymbol{n} + \dot{q}_{p,\mathrm{rad,net}}$ & — \\
        \hline
         $\partial\Omega_{g-\mathrm{room}}$ & — & $\lambda_{g}\nabla T_{g} \cdot \boldsymbol{n}=h_{\mathrm{ext}}(T_{\mathrm{room}}-T_g)$ & — \\
        \hline
         $\partial\Omega_{p-\mathrm{room}}$ & $\boldsymbol{u}_p=0$ & $\lambda_{p}\nabla T_{p} \cdot \boldsymbol{n}=0$ & — \\
        \hline
        \end{tabular}
    \caption{Boundary conditions for the glass–air–PCM coupled system.}
    \label{tab:boundaryc}
\end{table}

Initially, the full domain is assumed to be in thermal equilibrium at the ambient temperature, $T_{i0}=T_{\mathrm{amb}}=22\,^{\circ}$C. 
Accordingly, the air is at rest ($\boldsymbol{u}_{f0}=0$~m/s, $P_f=P_0$), and the PCM is completely solid ($\boldsymbol{u}_{p0}=0$~m/s, $P_p=0$).

\section{Numerical strategy}\label{sec:Numerical} 

\subsection{Forward full-order solution: finite volume method} 

The governing equations are solved using the finite volume method (FVM). 
The computational domain is discretized on a structured, anisotropic Cartesian mesh with a staggered-grid arrangement, where volume-averaged scalar quantities are stored at cell centers and velocity components at cell faces. 
Velocity and pressure coupling is achieved through the SIMPLEC pressure-correction algorithm, while nonlinearities are tackled via fixed-point iterations. 
The resulting algebraic system is solved with the bi-conjugate stabilized gradient (Bi-CGSTAB) method. 
Convergence is accepted when the maximum normalized residual of all state variables $\phi = (\boldsymbol{u}, T, P)$ drops below $10^{-6}$.

Convective fluxes of each transported variable at the cell faces, $\phi_S$, are reconstructed using the third-order weighted essentially non-oscillatory (WENO3) scheme, which ensures high-order accuracy in smooth regions and suppresses non-physical oscillations near sharp gradients. 
The WENO3 reconstruction combines two adaptive weights, $w_0$ and $w_1$ ($w_k \ge 0$, $w_0 + w_1 = 1$), and two-point sub-stencils $S_i^0(\phi_{i-1}, \phi_i)$ and $S_i^1(\phi_i, \phi_{i+1})$, according to
\[
\phi_{S_i} = w_0 S_i^0 + w_1 S_i^1,
\]
where the weights depend on local smoothness functions of $\phi$.
Diffusive fluxes at cell faces are computed using a second-order central-difference scheme for gradient reconstruction.

The temporal domain is divided into $N_t$ uniform intervals of size $\Delta t = t^{n+1} - t^n$. 
Time derivatives are discretized using an optimized backward differentiation formula (BDF2-opt), which is unconditionally stable and employs four time levels:
\begin{equation}
\partial_t \phi = \frac{c_1 \phi^{n+1} + c_2 \phi^n + c_3 \phi^{n-1} + c_4 \phi^{n-2}}{\Delta t}.
\end{equation}
The coefficients $c_1$–$c_4$ are obtained as a linear combination of the standard second- and third-order BDF schemes:
\[
c_i = 0.52\, c_i^{\text{BDF2}} + 0.48\, c_i^{\text{BDF3}}, \qquad i = 1, \ldots, 4.
\]

\subsubsection{Discretization study} 

The optimal spatial and temporal discretization was determined by evaluating the root-mean-square (RMS) error of velocity and temperature fields between four mesh sizes and four time-step sizes. 
Table~\ref{tab:discretization_study} summarizes the time-averaged RMS error between a given resolution $M$ and its immediate predecessor $M-1$, defined as
\[
RMS(\phi) = 
\left[ \frac{1}{N_v}  
\sum_{i=1}^{N_v} 
\left( \phi_{M}(i) - \phi_{M-1}(i) \right)^2 
\right]^{1/2},
\]
for both the air and PCM domains. 
In the spatial refinement study, the finer mesh was linearly interpolated onto the coarser one before computing the RMS values.

Each domain was independently meshed using a wall-concentrated hyperbolical tangent refinement, ensuring sufficient grid resolution to capture near-wall gradients and domain-specific flow dynamics. 
The time step was kept constant during most of the simulation, except for the initial 30~minutes, where it was reduced by half to better resolve the highly-dynamic behavior during the initial stage.

\begin{table}[h!]
\centering
    \begin{tabular}{cccccc}
    \toprule
     & Size & $RMS(\boldsymbol{u})$ in $\Omega_f$ & $RMS(T)$ in $\Omega_f$ & $RMS(\boldsymbol{u})$ in $\Omega_p$ & $RMS(T)$ in $\Omega_p$ \\
     & (minimum) & (m\,s$^{-1}\!\times\!10^{-3}$) & ($^{\circ}\text{C}\times\!10^{-3}$) & (m\,s$^{-1}\!\times\!10^{-3}$) & ($^{\circ} \text{C}\times10^{-3}$) \\
    \midrule
    $M_{2}$ & $9.1 \times 10^{-4}$ m & $8.10$ & $30.63$ & $3.32$ & $102.53$ \\
    $M_{3}$ & $3.8 \times 10^{-4}$ m & $5.77$ & $22.80$ & $1.95$ & $76.83$ \\
    $M_{4}$ & $1.4 \times 10^{-4}$ m & $5.18$ & $21.61$ & $1.66$ & $72.07$ \\
    \midrule
    $\Delta t_{2}$ & $6 \times 10^{-2}$ s & $2.29$ & $12.86$ & $0.66$ & $25.32$ \\
    $\Delta t_{3}$ & $4 \times 10^{-2}$ s & $1.56$ & $10.74$ & $0.43$ & $21.77$ \\
    $\Delta t_{4}$ & $3 \times 10^{-2}$ s & $1.41$ & $10.23$ & $0.38$ & $19.88$ \\
    \bottomrule
    \end{tabular}
    \caption{Comparison of RMS errors between successive mesh sizes ($M_1$–$M_4$) and time-step sizes ($\Delta t_1$–$\Delta t_4$) for the air domain ($\Omega_f$) and PCM domain ($\Omega_p$).}
    \label{tab:discretization_study}
\end{table}

When refining the mesh from $M_2$ to $M_3$, the RMS errors in velocity and temperature decrease by approximately 40\% and 25\%, respectively, whereas the transition from $M_3$ to $M_4$ yields only 12\% and 5\% improvements. 
These results indicate that further mesh refinement provides marginal accuracy gains. Therefore, mesh $M_3$ was selected for the remaining simulations.
Similarly, decreasing the time-step size from $\Delta t_2$ to $\Delta t_3$ reduces the velocity and temperature RMS errors by 32\% and 15\%, while the change from $\Delta t_3$ to $\Delta t_4$ leads to smaller variations of 11\% and 7\%. 
Considering that smaller time steps do not significantly affect the overall computational cost, $\Delta t_4$ was adopted for all subsequent simulations.

\subsubsection{Validation} 

\begin{figure}[h]
    \centering
    \includegraphics[width=1\linewidth]{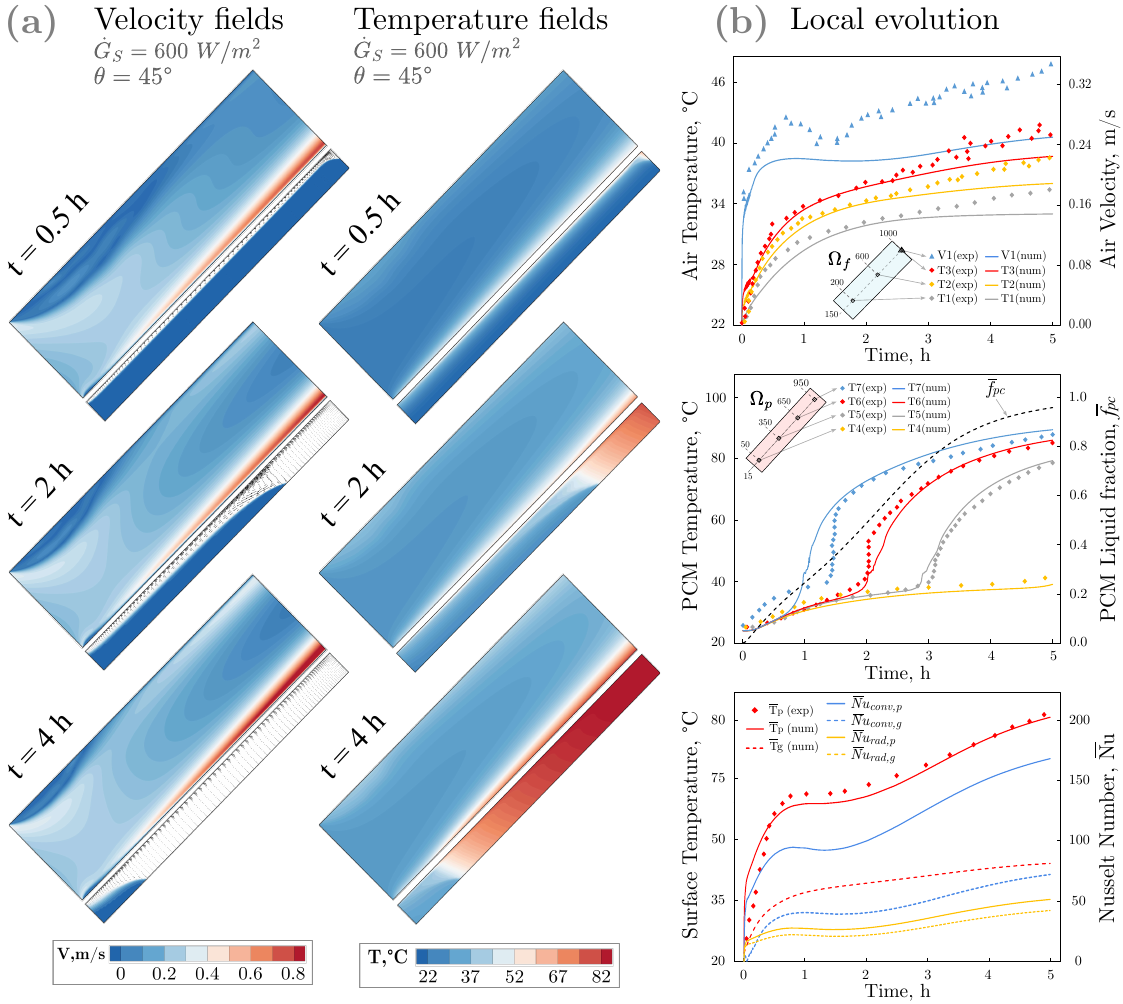}
    \caption{(a) Snapshots of velocity magnitude and temperature in the airflow and PCM domains ($\Omega_f$ and $\Omega_p$) at three different time instants. 
    (b) Local evolution of temperature and velocity at fixed points in $\Omega_f$ (top plot) and of temperature in $\Omega_p$ (middle plot), and of surface-averaged temperature and Nusselt numbers (bottom plot).
    Lines represent the forward numerical solution and symbols correspond to measurements from \cite{HUANG2024130154}.}
    \label{fig:forward_figures}
\end{figure}

This subsection analyzes the main features of the forward full-order solution and assesses its physical fidelity against the available experimental data. 
Figure~\ref{fig:forward_figures}(a) presents representative snapshots of velocity and temperature fields within both the airflow and PCM domains ($\Omega_f$ and $\Omega_p$), while Figure~\ref{fig:forward_figures}(b) compares the temporal evolution of local temperatures at selected points and the surface-averaged quantities obtained from simulations and experiments.

In the air domain, the velocity field exhibits two characteristic flow features. 
First, the inflow enters the cavity with a marked inclination toward its center, generating a large recirculation region along the upper (glass) wall and a smaller one along the lower (PCM) wall. 
These vortical structures shed periodically, producing oscillations that dissipate downstream along the air channel. 
During the first hour, the upper recirculation zone diminishes rapidly as buoyancy forces accelerate the near-wall flow, thereby reducing both the size and intensity of the vortices. 
Second, a clear asymmetry develops between the upper and lower boundary layers, caused not only by the collector's inclination but also by the thermal boundary conditions: the lower wall directly absorbs the incident solar flux $\dot{G}_S$, while the upper wall receives only a fraction of this energy via surface radiation. 
The temperature field mirrors these effects: high-temperature isotherms detach toward the cavity core due to vortex shedding, and the lower wall exhibits stronger stratification of warm isotherms than the upper wall.

In the PCM domain, melting initiates near the upper region and progresses predominantly from the right-hand wall as the charging stage advances. 
Velocity vectors intensify in the intermediate zone, where internal convection becomes more pronounced. 
Isotherms display a moderate inclination and propagate steadily from the right boundary, except in the central region where recirculating flow distorts their distribution. 
This behavior arises from buoyancy-driven internal motion, which advects the high-temperature region toward the upper-right corner of the PCM.

When compared with experimental measurements, the forward simulation reproduces the main spatio-temporal trends of temperature and velocity but underestimates their magnitude, particularly near the outlet. 
The largest deviations occur toward the end of the charging period ($t>3$~h), where the experimental data exhibit strong fluctuations attributed to turbulent effects resulting from boundary-layer separation. 
In the PCM region, numerical predictions show better agreement during the final phase of melting, although discrepancies remain near the inlet and in the early transition from solid to liquid. 
Despite these differences, the surface-averaged temperature of the PCM plate agrees well with the experimental trends. 
Finally, the computed convective and radiative Nusselt numbers, normalized by $\dot{q}_{\mathrm{cond}} = \lambda_f \Delta T_f / H$, indicate that radiation contributes as significantly as convection to the overall heat-transfer process, particularly on the glass surface during the initial stage ($t<0.5$~h).

\subsection{Data assimilation framework} 

The workflow consists of four main stages:
(1) a reduced-order modeling stage that generates a time-dependent, parametric reduced basis from a dataset of forward snapshots;
(2) a data-filling stage in which the initial set of measurements is expanded into two augmented datasets using both data- and physics-driven interpolation methods; 
(3) a least-squares reconstruction stage for dynamic temperature fields using real and synthetic measurements; 
and (4) a velocity field estimation stage that incorporates data-corrected buoyancy-source terms to improve the prediction of the outlet airflow velocity.

\subsubsection{Least-squares reconstruction with ROM regularization} 

Let $\ell\in\mathbb{R}^{m}$ denote the vector of measurements and let $W\in\mathbb{R}^{\mathcal{N}\times m}$ represent the linear observation operator such that the full-order model predictions at the measurement locations are given by $W^{\top}T\in\mathbb{R}^{m}$ for a state vector $T\in\mathbb{R}^{\mathcal{N}}$. 
Because typically the number of measurements is much lower than the number of full-order states, $m\ll \mathcal{N}$, the inverse problem is not well-posed. To regularize it, the solution is constrained to a subspace spanned by the reduced basis $\Phi \in \mathbb{R}^{\mathcal{N}\times n}$, with $n \ll \mathcal{N}$, so that the state can be expressed as $T = \Phi c$. Hence, the regularized least-squares optimization problem is
\begin{equation}
    c^\star=\arg\min_{c\in\mathbb{R}^{n}}\ \big\|\ell-W^{\top}\Phi c\big\|_2^2.
\end{equation}
This formulation leads to the normal equations, $(G^{\top}G)c=G^{\top}\ell$, where $G := W^{\top}\Phi \in \mathbb{R}^{m\times n}$ is a cross-Gramian matrix, and the reconstructed state field is $T^* = \Phi c^\star$. This subspace serves as a physics-informed regularization and admits an efficient offline/online split: $\Phi$ and $W$ are assembled offline, while in the online stage, only a small $n\times n$ system is solved online (\cite{galarce_2025_ijmecsci,Maday_2015_ESAIM}).

The reduced-order basis $\Phi$ is constructed via a truncated singular value decomposition (SVD) to a snapshot matrix, $A=U_n\Sigma_n V_n^{\top}$, where the number of modes ($n$) defines the projection space capturing the dominant dynamics of the thermo-fluid system required for accurate state estimation (\cite{Willcox_2006_jcompfluid}). The matrix $A$ is assembled from a dataset of forward simulations spanning a parameter space $\theta_k$, such that: $A=[\,T(\theta_1)\ |\ \cdots\ |\ T(\theta_K)\,]\in\mathbb{R}^{\mathcal{N}\times K}$. The procedure used to generate the forward snapshot dataset and determine the optimal number of modes is detailed in Section \ref{sec:Results_ROM}.

An \textit{a-priori} error bound can be derived before performing the least-squares reconstruction, which enables determining the optimal number of modes in the reduced basis for a given sensor configuration (\cite{Maday_2015, galarce2023_MRE}):
\begin{equation} 
    \frac{||T-T^*||}{T} \le \frac{E_{\Phi}}{\hat{S}_n} = e_p(n), \quad E_{\Phi} =\frac{||A-A_n||_{FRO}}{||A||_{FRO}},\quad \hat{S}_n = \inf_{c \in \mathds{R}^n} \frac{\| \boldsymbol{G} c\|}{\| \Phi c\|}
\end{equation}
where $E_{\Phi}$ represents the POD truncation error as the reduced basis's size increases, and $\hat{S}_n$ is the last singular value of the cross-Gramian matrix $\boldsymbol{G}$, which quantifies the observability of the sensor configuration. This error bound indicates a trade-off between the increasing quality of the ROM and the loss of observability as the size of the reduced basis increases, a critical issue for well-posedness in data assimilation problems.

\subsubsection{Data filling procedure} 

\begin{figure}[h]
    \centering
    \includegraphics[width=1\linewidth]{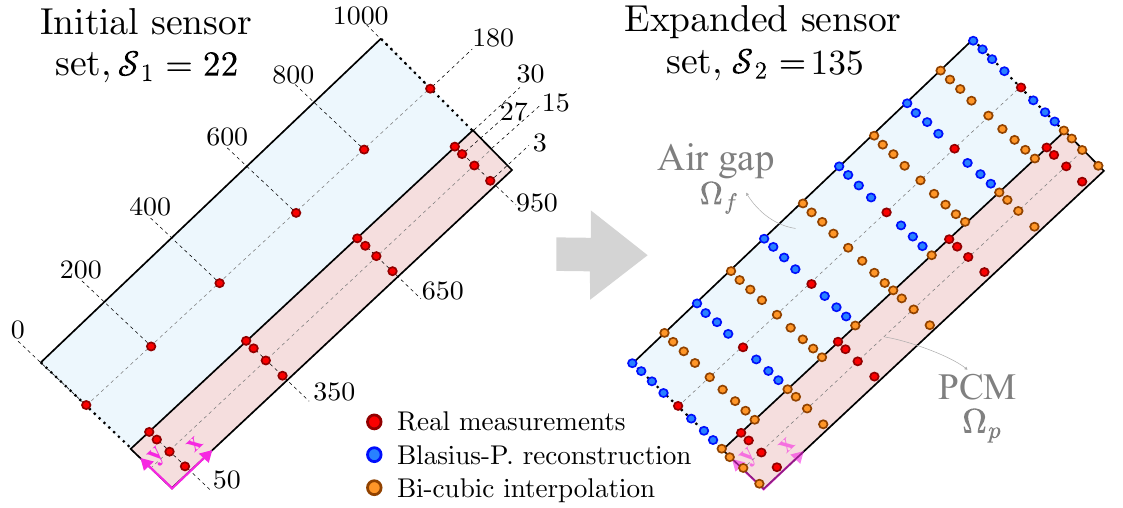}
    \caption{Positions of the initial measurement set $\mathcal{S}_1$ (left) and the first expanded set $\mathcal{S}_2$ (right) in both the air and PCM domains.}
    \label{fig:setofsensors}
\end{figure}

The initial measurement set, $\mathcal{S}_1 = 22$ points (six located in the air domain $\Omega_f$ and sixteen in the PCM domain $\Omega_p$) is expanded to $\mathcal{S}_2 = 135$ and $\mathcal{S}_3 = 203$ total points using two complementary interpolation techniques that combine data-driven and physics-based approaches, illustrated in Figure~\ref{fig:setofsensors}. 
A brief description of each method follows:
\begin{itemize}
    \item \textbf{Bi-cubic interpolation:} 
    A data-driven spline operator used to reconstruct a smooth temperature field over a uniform grid via the tensor product of one-dimensional cubic splines in $x$ and $y$. 
    The spline coefficients are obtained by solving tridiagonal systems that ensure $C^2$ continuity along grid lines. Within each cell, the temperature distribution is approximated by a bi-cubic polynomial. 
    For a grid spacing $h = x_{i+1} - x_i$ and second derivatives $y''_i$, the local spline on $[x_i, x_{i+1}]$ is defined as:
    \[
    s(x) = a\,y_i + b\,y_{i+1} + \frac{h^2}{6} \Big[ (a^3 - a)\,y''_i + (b^3 - b)\,y''_{i+1} \Big],
    \]
    where $a$ and $b$ are normalized local coordinates satisfying $a+b=1$.
    
    \item \textbf{Blasius–Pohlhausen similarity solution:} 
    A physics-driven technique that reconstructs laminar thermal boundary-layer profiles representing forced convection over a flat plate.  The governing equations are reduced to ordinary differential equations (ODEs) using the similarity variable $\eta = y \sqrt{U_\infty / \nu x}$. 
    The velocity profile satisfies Blasius’s equation,
    \[
    f''' + \tfrac{1}{2} f f'' = 0,
    \]
    while the normalized temperature follows
    \[
    \theta'' + \tfrac{1}{2} Pr\, f(\eta)\, \theta' = 0, \qquad \theta(0)=1, \quad \theta(\infty)=0,
    \]
    where $\theta(\eta) = (T - T_\infty)/(T_w - T_\infty)$. 
    The resulting ODE system is integrated using a fourth-order Runge–Kutta (RK4) scheme, and the dimensional temperature field is reconstructed as 
    $T(y) = T_\infty + (T_w - T_\infty)\, \theta(\eta)$.
\end{itemize}

The construction of the expanded measurement sets $\mathcal{S}_2$ and $\mathcal{S}_3$ is carried out as follows:  
(i) first, the bi-cubic interpolation is applied in $\Omega_p$, where the initial measurement density is higher;  
(ii) subsequently, the Blasius–Pohlhausen solution is used to reconstruct temperature profiles in $\Omega_f$ across the glass and PCM boundary layers, where $U_\infty$ corresponds to the measured outlet velocity of the air gap; and (iii) the bi-cubic interpolation is applied again to interpolate intermediate regions in $\Omega_f$.

\section{Results and discussion}\label{sec:Results} 

This section describes the construction of the reduced-order basis and evaluates its POD truncation error. We also compute the corresponding \textit{a-priori} error bounds for each dataset. Subsequently, the ROM-based data assimilation framework is assessed using synthetic temperature measurements generated from the full-order solution. Finally, we apply the ROM-DA framework to the experimental datasets to reconstruct dynamic temperature fields, which are then assimilated into an enhanced velocity-field estimation. The data assimilation framework and the construction of the ROM was implemented using the in-house MAD library \cite{galarceThesis, MELLA2026108145}.

\subsection{Reduced-order modeling} \label{sec:Results_ROM} 
\begin{figure} [h]
    \centering
    \includegraphics[width=1\linewidth]{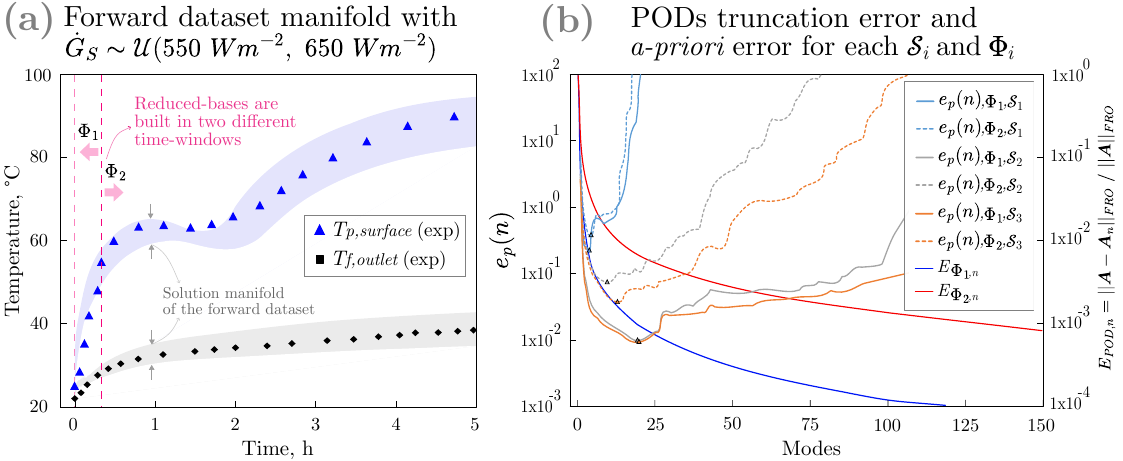}
    \caption{(a) Solution manifold of the forward dataset. Symbols are experimental data in fixed points at the surface of the PCM and at the outlet air gap, while the shaded regions represent maximum and minimum values of the forward dataset. (b) POD truncation error decay of each reduced basis ($\Phi_1$ and $\Phi_2$), and \textit{a-priori} error bounds for each reduced basis and set of sensors.}
    \label{fig:ROM_plots}
\end{figure}

The full-order snapshot matrix $A$ is assembled from multiple forward simulations with parameters randomly sampled in the vicinity of the experimental conditions. A total of 100 forward simulations were carried out, each producing 960 snapshots uniformly distributed in time, with solar radiation varied with a random uniform distribution within $\pm 12\%$ around the nominal value, $\dot{G}_S \sim \mathcal{U}(550\ \mathrm{Wm}^{-2},\ 650\ \mathrm{Wm}^{-2})$. 
This sampling strategy ensures that the generated solution manifold adequately encompasses the experimental observations, as shown in  Figure~\ref{fig:ROM_plots}(a).

To improve the representation of the underlying physics of the process, two reduced bases $\Phi_1$ and $\Phi_2$ (or POD) are constructed independently according to two characteristic time-windows, in which distinct dominant time scales are identified: (1) for $t < 0.3$~h, where the state variables increase rapidly from the initial condition and significant fluctuations arise due to vortex shedding; and (2) for $t > 0.3$~h, where the variables remain more stable over time and fluctuations are negligible given the sampling resolution. By doing this, we achieve reduced bases with truncation errors that rapidly decay to values below 10$^{-2}$ with fewer than 25 modes, as shown in Figure~\ref{fig:ROM_plots}(b).

Figure~\ref{fig:ROM_plots}(b) also presents the \textit{a-priori} error bounds computed for each reduced basis and sensor set. For $\Phi_1$, the optimal numbers of modes are 4, 19, and 19 for $\mathcal{S}_1$, $\mathcal{S}_2$, and $\mathcal{S}_3$, respectively, whereas for $\Phi_2$, the corresponding values are 4, 10, and 14.
These results confirm that the observability and truncation errors exhibit consistent trends across both reduced bases and measurement sets.

\subsection{Data assimilation with synthetic measurements} \label{sec:Results_DA_synthetic} 

To assess the reconstruction accuracy of the ROM-DA framework, we first compare the dynamic temperature field reconstruction using synthetic measurements. These measurements are generated from a set of 16 forward full-order simulations that are independent of those used to construct the reduced bases $\Phi_1$ and $\Phi_2$. The forward solution of each simulation is sampled at discrete sensor locations with a diameter of 0.01~m, matching the positions of the three measurement sets ($\mathcal{S}_1$, $\mathcal{S}_2$, and $\mathcal{S}_3$). In this way, the synthetic measurements replicate the limited spatial resolution and coverage of the experimental data while remaining free from measurement noise or bias typically associated with instrumentation.

\begin{figure}[h]
    \centering
    \includegraphics[width=1\linewidth]{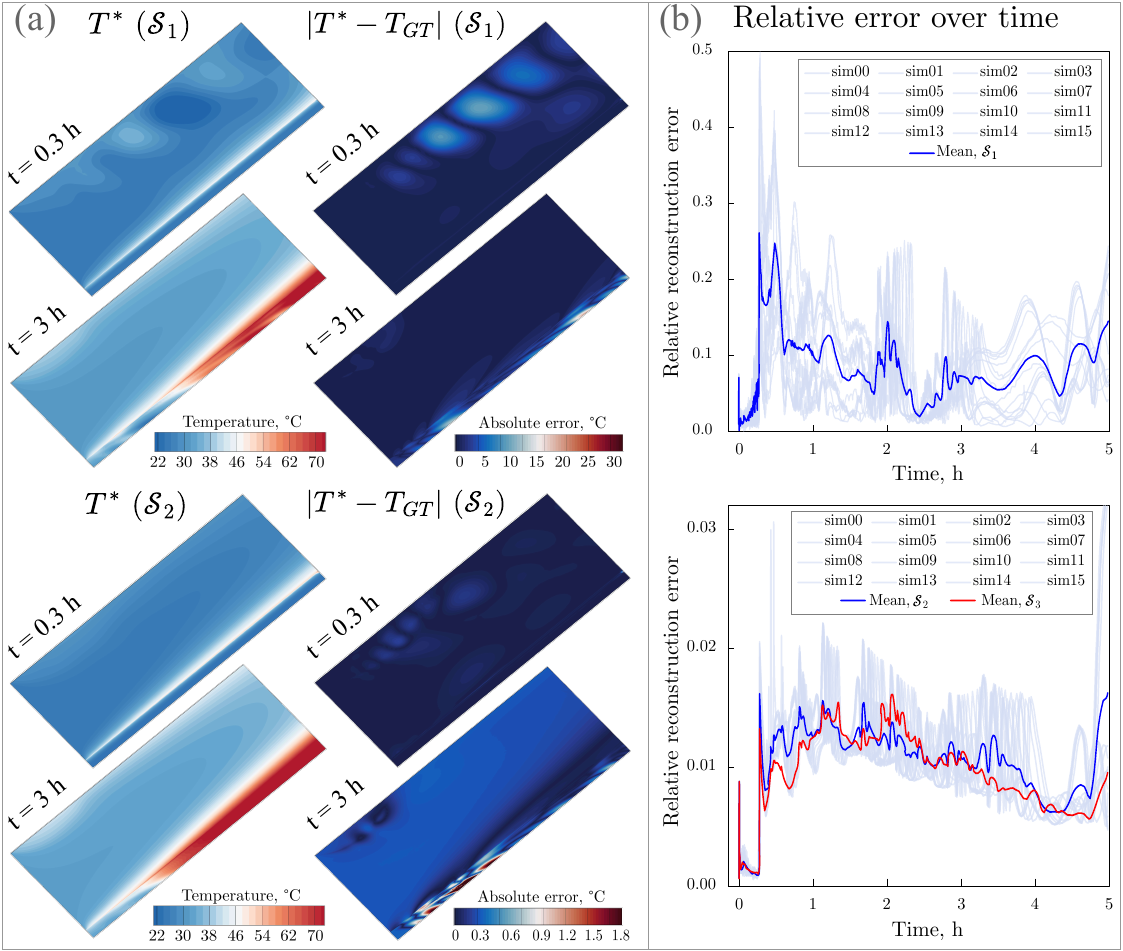}
    \caption{(a) Snapshots of the reconstructed temperature fields $T^*$ and of their absolute errors relative to the full-order solution (or ground through, $T_{GT}$), using sensor sets $\mathcal{S}_1$ and $\mathcal{S}_2$.(b) Relative reconstruction error over time for the test dataset comprising the 16 forward simulations and their mean error over time for each sensor set.}
    \label{fig:synthetic_reco}
\end{figure}

Figure~\ref{fig:synthetic_reco}(a) compares two representative snapshots of the reconstructed temperature fields, corresponding to the simulation with the highest mean relative error, together with their absolute deviations from the full-order (ground-truth) solution, $T_{GT}$, for sensor sets $\mathcal{S}_1$ and $\mathcal{S}_2$. The temperature field reconstructed with $\mathcal{S}_1$ shows pronounced local errors within the air domain ($\Omega_f$) at the initial time instant ($t = 0.3$~h), with deviations up to 15~\textdegree C near the upper-wall vortex-shedding region. In contrast, at the later time instant ($t = 3$~h), the errors are mainly concentrated in the PCM domain ($\Omega_p$), reaching values above 30~\textdegree C close to the phase-change fronts. For the reconstruction obtained with $\mathcal{S}_2$, the absolute errors are substantially lower, yielding smoother temperature fields that closely match the forward solution. During the first stage, the dominant discrepancies occur in $\Omega_f$ near the vortex-shedding zone, remaining below 1~\textdegree C, whereas during the second stage they appear in $\Omega_p$ with magnitudes up to 2~\textdegree C. This contrast in error magnitude between the two sensor sets indicates that $\mathcal{S}_1$ provides insufficient spatial information, causing the least-squares minimization to generate higher-amplitude oscillations in regions with strong nonlinear dynamics.

To quantify these errors, Figure~\ref{fig:synthetic_reco}(b) shows the full-domain-averaged relative reconstruction error over the charging process for the test dataset comprising 16 forward simulations, together with the mean errors associated with each sensor set. The dynamic trends confirm the previously discussed differences: with $\mathcal{S}_1$, the time-averaged relative error reaches approximately 10\%, whereas for $\mathcal{S}_2$ and $\mathcal{S}_3$, the values drop to around 1\%. Moreover, two consistent tendencies emerge: (i) the relative errors within the interval corresponding to the first reduced basis $\Phi_1$ are considerably smaller than those in the second interval, and (ii) increasing the number of sensors from $\mathcal{S}_2$ to $\mathcal{S}_3$ yields only marginal improvements. The first observation confirms that the temporal separation of reduced bases enhances reconstruction accuracy during the early transient, while the second demonstrates that beyond a certain measurement density, additional sensors provide diminishing returns in reconstruction performance.

\subsection{Data assimilation with real measurements} \label{sec:Results_DA_real} 

The previous section demonstrated that the proposed ROM-DA framework accurately reconstructs dynamic temperature fields using synthetic measurements outside the ROM training dataset. 
When applied to real measurements, however, additional uncertainty sources arise that are absent in idealized numerical tests. 
While synthetic data generated under white Gaussian noise preserves full consistency between the model, observation operator, and noise statistics, real measurements are affected by instrumental noise, sensor bias, temporal drift, and spatial-averaging errors, and the physical system itself may deviate from the model due to parameter or model-form discrepancies~\cite{GALARCE_2025_appliedmathmodelling}. 
In this section, we assess the performance of the framework under these realistic conditions, without applying explicit noise or bias corrections. 
For each reduced basis and sensor set, the number of modes was selected to maximize agreement with experimental data; although these values differ slightly from those predicted by the \textit{a-priori} error analysis, they follow the same overall trend.

\begin{figure}[h]
    \centering
    \includegraphics[width=1\linewidth]{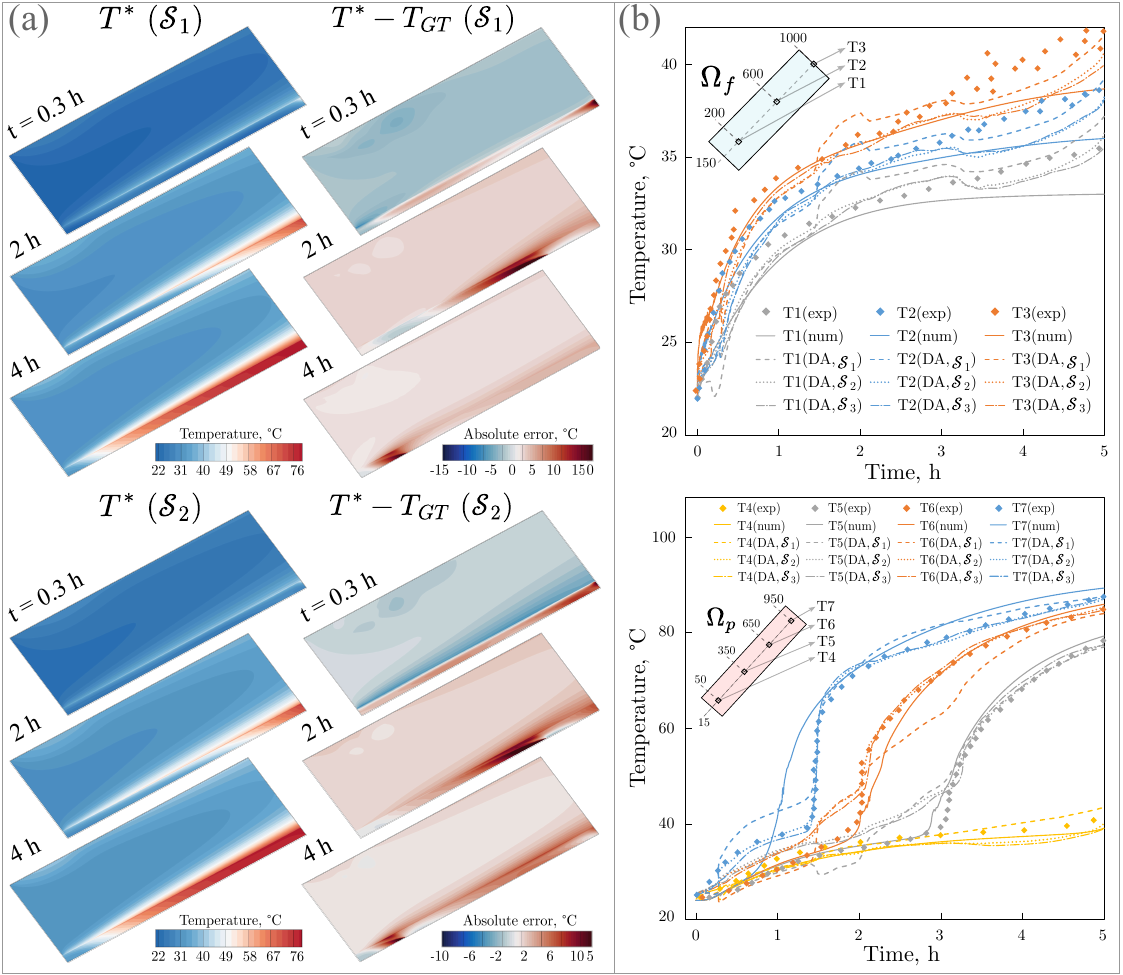}
    \caption{(a) Snapshots of reconstructed temperature field $T^*$ and relative error with respect to the forward solution (ground truth, $T_{GT}$) for sensor sets $\mathcal{S}_1$ and $\mathcal{S}_2$. (b) Comparison of local temperature evolution at fixed points for the experimental data (symbols), the forward model (continuous line), and for the reconstructed with each sensor set (dashed lines).}
    \label{fig:real_reco}
\end{figure}

Figure~\ref{fig:real_reco}(a) compares three representative snapshots of the reconstructed temperature field, $T^*$, and their corresponding relative error fields with respect to the forward solution, $T_{GT}$, for sensor sets $\mathcal{S}_1$ and $\mathcal{S}_2$. The reconstructed temperature fields show no visually distinguishable differences between both sets, but the relative error fields reveal recurring spatial patterns. 
At the initial instant ($t = 0.3$~h), $T^*$ is lower than the forward solution across the entire air domain ($\Omega_f$) but higher in the PCM domain ($\Omega_p$). At later instants, $T^*$ in $\Omega_f$ becomes higher than the forward solution, particularly near the wall boundary layers, reaching deviations up to 6~\textdegree C for $\mathcal{S}_1$ and 4~\textdegree C for $\mathcal{S}_2$. In $\Omega_p$, the fluctuations are larger than in $\Omega_f$, with errors ranging from –6~\textdegree C to 18~\textdegree C for $\mathcal{S}_1$ and from –4~\textdegree C to 12~\textdegree C for $\mathcal{S}_2$. The differences between $T^*$ at successive instants indicate that the temporal separation of reduced bases influences each time window differently, slightly affecting the continuity of the reconstruction. 
Moreover, the alternating positive and negative deviations within the PCM domain suggest greater uncertainty during the liquid–solid transition, where distinct physical regimes coexist.

Figure~\ref{fig:real_reco}(b) compares the temporal evolution of $T^*$ at selected locations against both experimental and forward-simulation data. 
Overall, all sensor sets improve the agreement between the numerical predictions and the experimental measurements. However, all three show a clear discontinuity at $t = 0.3$~h, corresponding to the switching point between the two reduced bases, and inconsistencies in the outlet region—where measurement noise is more significant—as well as during the phase-change intervals.

In the air domain ($\Omega_f$), the reconstruction obtained with $\mathcal{S}_1$ displays larger deviations and lower spatial uniformity during the first interval ($t < 0.3$~h), yet it predicts higher outlet temperatures in the second interval, achieving better fidelity than $\mathcal{S}_2$ and $\mathcal{S}_3$. This outcome is somewhat unexpected, since—as shown in Section~\ref{sec:Results_DA_synthetic}—a higher sensor density generally produces reconstructions with lower relative errors. 
Additionally, increasing the number of sensors from $\mathcal{S}_2$ to $\mathcal{S}_3$ provides no noticeable improvement in overall accuracy. 
Conversely, in the PCM domain ($\Omega_p$), sets $\mathcal{S}_2$ and $\mathcal{S}_3$ substantially improve the temperature predictions near the inlet, where discrepancies were previously larger, while slightly worsening them near the outlet region.

Since the primary function of a solar chimney is to induce a buoyancy-driven ventilation flow, the main performance indicator is the volumetric flow rate, which can be directly estimated from the outlet velocity. Therefore, the dynamic temperature field, although informative, is not key design quantity. Here, we assess whether the temperature-corrected fields can improve the prediction of the velocity field obtained from the thermally coupled Navier–Stokes equations. Because the flow is driven exclusively by temperature gradients, the reconstructed fields are assimilated through the buoyancy term. Given that the temporal resolution of $T^*$ is lower than the needed by the forward solver, a cubic temporal interpolation is applied to obtain the instantaneous buoyancy force.

\begin{figure}[h]
    \centering
    \includegraphics[width=1\linewidth]{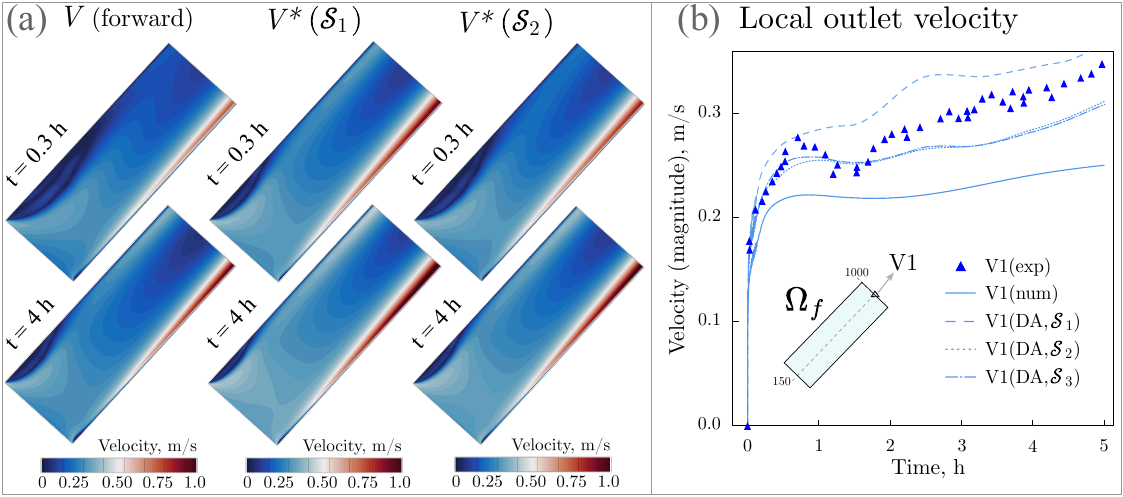}
    \caption{(a) Velocity snapshots of the forward model and the data-corrected method with $\mathcal{S}_1$ and $\mathcal{S}_2$. (b) Comparison of local outlet velocity between the different models and experimental data.}
    \label{fig:velocity_reco}
\end{figure}

Figure~\ref{fig:velocity_reco}(a) compares the velocity snapshots from the forward model with those computed using the data-corrected temperature fields at two characteristic instants for $\mathcal{S}_1$ and $\mathcal{S}_2$. 
The overall flow structure remains similar to that of the forward model, with the main differences being a reduction in vortex size near the inlet corners and a thicker boundary layer along the upper wall. However, the magnitude of the velocity field increases, particularly for $\mathcal{S}_1$, where the peak values exceed 1~m/s within the boundary layer at the initial instant. 
These results are consistent with the reconstructed temperatures fields with $\mathcal{S}_1$ that are notably higher than the forward solution, which increases the buoyancy forces and thus accelerate the airflow.

Figure~\ref{fig:velocity_reco}(b) compares the temporal evolution of the local outlet velocity—computed as the average over a 0.05~m-radius region centered at the midplane, matching the spatial resolution of the experimental sensor—between the experimental data, the forward model, and the data-corrected predictions for all sensor sets. 
The dynamic velocity predictions obtained with the assimilated temperature fields exhibit substantial improvement compared with the forward model: the time-average root-mean-square (RMS) error with the experimental data decreases by 44\% for $\mathcal{S}_1$ and by 55\% for $\mathcal{S}_2$ and $\mathcal{S}_3$. 
Moreover, all three sensor sets better capture the transient velocity variations observed between 1~h~$\le t \le$~2~h and for $t > 2$~h. 
Consistent with the thermal analysis, the $\mathcal{S}_1$-based reconstruction slightly overestimates the outlet velocity by approximately 18\% relative to the experiments, whereas $\mathcal{S}_2$ and $\mathcal{S}_3$ yield modest underestimations of about 12\%.

\section{Conclusions} \label{sec:Conclusions} 

This study assessed a reduced-order model data-assimilation (ROM–DA) framework to improve the prediction accuracy of outlet airflow velocity in an inclined PCM-integrated solar chimney. The framework combines a reduced-order basis generated from a full-order solver (CFD) with three measurement sets expanded through a hybrid data-filling strategy. The reconstructed temperature fields are then assimilated into the forward solver to obtain data-corrected dynamic velocity fields.

The ROM–DA framework was successfully validated with synthetic measurements derived from the full-order model, accurately reconstructing the dynamic temperature fields in both the air and PCM domains. The maximum relative errors remained below 10\% for $\mathcal{S}_1$ and below 3\% for $\mathcal{S}_2$ and $\mathcal{S}_3$. Minor deviations observed near vortex-shedding regions and phase-change fronts were attributed to the dynamic basis switching and the inherent challenges of multi-domain reconstruction.

When assimilating real experimental data, the ROM–DA framework enhanced the physical fidelity of the temperature evolution relative to the forward solver. It captured transient variations during the phase-change process and toward the end of the heating stage, although discrepancies persisted in regions affected by measurement noise or strong nonlinearities. Among the sensor configurations, $\mathcal{S}_1$ produced less uniform reconstructions but yielded slightly higher temperatures near the outlet, resulting in improved local temperature predictions compared to $\mathcal{S}_2$ and $\mathcal{S}_3$.

The thermally coupled flow solver, augmented with assimilated temperature fields, improved the prediction of outlet velocity without requiring direct velocity measurements. By correcting buoyancy forces using assimilated data, the framework reduced the RMS error of outlet velocity by 44\% for $\mathcal{S}_1$ and by 55\% for $\mathcal{S}_2$ and $\mathcal{S}_3$ relative to the forward model. The higher reconstructed temperatures from $\mathcal{S}_1$ slightly intensified buoyancy effects, leading to a mild overprediction of velocity magnitude.

These results demonstrate that the proposed ROM–DA framework provides a reliable and computationally efficient strategy for reconstructing primary states with better physical fidelity than traditional methods in multi-domain thermo-fluid systems. The hybrid data-filling approach effectively supports the assimilation process, enabling accurate predictions even under sparse or noisy measurements. 
The practical implications of this work are: (i) generation of a numerical tool that allows optimizing and evaluating solar collector devices under different operating or environmental conditions, (ii) improving the design criteria of sensor networks for efficient monitoring of thermo-fluid systems. Future work will focus on two aspects: (i) generation of enhanced reduced bases both in dynamic time-windows and in multi-domains with different governing physics; (ii) coupling this methodology with local meteorological measurements to generate ventilation efficiency maps that can be updated in near-real-time.

\hfill

\textit{Acknowledgments}

E.~Castillo and F.~Galarce acknowledge the financial support of the Chilean National Agency for Research and Development (ANID) through the project Fondecyt No. 1250287. F.~Galarce also acknowledges DI Vinci PUCV Iniciación 039.731/2025 and the Horizon Europe - 2nd Opportunity OPPTY-MSCA/0125 funding from the Cyprus Research and Innovation Foundation. D.~Rivera acknowledges the Vice-Rectory for Research, Innovation and Creation (VRIIC) of Universidad de Santiago de Chile for the support in the Postdoctoral Proyect DICYT-POSTDOC 052416CD. D.R.Q.~Pacheco acknowledges funding by the Federal Ministry of Education and Research (BMBF) and the Ministry of Culture and Science of the German State of North Rhine-Westphalia (MKW) under the Excellence Strategy of the Federal Government and the Länder \\

\textit{Disclosure statement:} \\No potential conflict of interest is reported by the authors



\end{document}